\documentclass[fleqn]{mat01}
\usepackage{times,mathtimy,amssymb,latexsym}
\begin{document}

\newtheorem{theo}[defin]{\bf Theorem}
\newtheorem{propo}[defin]{\rm PROPOSITION}
\newtheorem{deff}[defin]{\rm DEFINITION}
\newtheorem{lem}[defin]{Lemma}
\newtheorem{rem}[defin]{\it Remark}
\newtheorem{coro}[defin]{\rm COROLLARY}

\markboth{Amir Khosravi and Behrooz Khosravi}{Frames and bases in
tensor products}

\title{Frames and bases in tensor products of Hilbert spaces\\ and Hilbert
$\pmb{C}^*$-modules}

\author{AMIR KHOSRAVI and BEHROOZ KHOSRAVI$^{*}$}


\address{Faculty of Mathematical Sciences and Computer
Engineering, University for Teacher Education, 599 Taleghani Ave.,
Tehran 15614, Iran\\
\noindent $^{*}$Department of Pure Mathematics, Faculty  of
Mathematics and Computer Science, Amirkabir University of
Technology (Tehran
Polytechnic), 424, Hafez Ave., Tehran~15914, Iran\\
\noindent E-mail: khosravi$_{-}$amir@yahoo.com}

\volume{117}

\mon{February}

\parts{1}

\pubyear{2007}

\Date{MS received 17 November 2005; revised 9 December 2005}

\begin{abstract}
In this article, we study tensor product of Hilbert $C^*$-modules
and Hilbert spaces. We show that if $E$ is a Hilbert $A$-module
and $F$ is a Hilbert $B$-module, then tensor product of frames
(orthonormal bases) for $E$ and $F$ produce frames (orthonormal
bases) for Hilbert $A\otimes B$-module $E\otimes F$, and we get
more results.

For Hilbert spaces $H$ and $K$, we study tensor product of frames
of subspaces for $H$ and $K$, tensor product of resolutions of the
identities of $H$ and $K$, and tensor product of frame
representations for $H$ and $K$.
\end{abstract}


\maketitle

\section{Introduction}

Gabor \cite{Ga}, in 1946 introduced  a technique for signal
processing which eventually led to wavelet theory. Later in 1952,
Duffin and Schaeffer \cite{sh} in  the  context of nonharmonic
Fourier series introduced frame theory for Hilbert spaces. In
1986, Daubechies, \hbox{Grassman} and Meyer \cite{dgm}  showed
that Duffin and Schaeffer's definition was an abstraction of
Gabor's concept. Frames are used in signal processing, image
processing, data compression, sampling theory, migrating the
effect of losses in packet-based communication systems and
improving the robustness of data transmission. Since tensor
product is useful in the approxi\-mation of multi-variate
functions of combinations of univariate ones, Khosravi and Asgari
\cite{ka} introduced frames in tensor product of Hilbert spaces.
Meanwhile, the notion of frames in Hilbert $C^*$-modules was
introduced and some of their properties were investigated
\cite{fornasier,FL1,FL2,Heilw,Lance}. In this article, we study
the frames and bases in tensor product of Hilbert $C^*$-modules
which were introduced in \cite{Lance} and we generalize the
techniques of \cite{ka} to $C^*$-modules.

In \S 2, we briefly recall the definitions and basic properties of
Hilbert $C^*$-modules. In \S 3, we investigate tensor product of
Hilbert $C^*$-modules, which is introduced in \cite{Lance} and we
show that tensor product of frames  for Hilbert $C^*$-modules $E$
and $F$, present frames for $E\otimes F$, and tensor product of
their frame operators is the frame operator of the tensor product
of frames. We also show that tensor product of frames of subspaces
produce a frame of subspaces for their tensor product. In \S 4, we
study resolution of the identity and prove that tensor product of
any  resolutions of $H$ and $K$, is a resolution of the identity
for $H\otimes K$. In \S 5, we study the frame representation and
we show that tensor product of frame vectors is a frame vector.
Also we show that tensor product of analysis operators (resp.
decomposition operators) is an analysis operator (resp. a
decomposition operator).

Throughout this paper, $\Bbb N$ and $\Bbb C$ will denote the set
of natural numbers and the set of complex numbers, respectively.
$A$ and $B$ will be unital $C^*$-algebras.

\section{Preliminaries}

Let $I$ and $J$ be countable index sets. In this section we
briefly recall the definitions and basic properties of Hilbert
$C^*$-modules and frames in Hilbert $C^*$-modules.  For
information about frames in Hilbert spaces we refer to
\cite{cas.art,Heilw,chrote,young}. Our
 reference for $C^*$-algebras is \cite{mur,Weg-alg}. For a $C^*$-algebra $A$
if $a\in A$ is positive we  write $a\geq 0$  and $A^+$ denotes the
set of positive elements of $A$.

\begin{deff}$\left.\right.$\vspace{.5pc}

\noindent {\rm Let $A$ be a unital $C^*$-algebra and let $H$ be a
left $A$-module, such that the linear structures of $A$ and $H$
are compatible. $H$ is a {\it pre-Hilbert $A$-module} if $H$ is
equipped with an $A$-valued inner product
$\langle.,.\rangle\hbox{\rm :}\ H\times H \rightarrow A$, that is
sesquilinear, positive definite and respects the module action. In
other words,

\begin{enumerate}
\renewcommand{\labelenumi}{(\roman{enumi})}
\leftskip .5pc
\item $\langle x,x\rangle  \geq 0$ for all $x\in H$ and $\langle
x,x\rangle  =0$ if and only if $x=0$;

\item $\langle  ax+y,z\rangle  =a\langle  x,z\rangle  +\langle
y,z\rangle  $ for all $a\in A$ and $x,y,z\in H$;

\item $\langle  x,y\rangle  = \langle  y,x\rangle  ^*$ for all
$x,y\in H$.\vspace{-.3pc}
\end{enumerate}

For $x\in H$, we define $\|x\|=\|\langle  x,x\rangle  \|^{1/2}$.
If $H$ is complete with $\|.\|$, it is called a {\it Hilbert
$A$-module} or a {\it Hilbert $C^*$-module} over $A$. For every
$a$ in $C^*$-algebra $A$, we have $|a|=(a^*a)^{1/2}$  and the
$A$-valued norm on $H$ is defined by $|x|=\langle  x,x\rangle
^{1/2}$ for $x\in H$.}
\end{deff}

\begin{deff}$\left.\right.$\vspace{.5pc}

\noindent {\rm Let $H$ be a Hilbert $A$-module. A family
$\{x_i\}_{i\in I}$ of elements of $H$ is a {\it frame} for $H$, if
there exist constants $0<A\leq B<\infty$, such that for all $x\in
H$,
\begin{equation}
A\langle  x,x\rangle  \leq \sum_{i\in I}  \langle  x,x_i\rangle
\langle  x_i,x\rangle  \leq B\langle  x,x\rangle.
\end{equation}
The numbers $A$  and $B$ are called lower  and upper bound of the
frame, respectively. If $A=B=\lambda$, the frame is {\it
$\lambda$-tight}. If $A=B=1$, it is called a {\it normalized tight
frame} or a {\it Parseval} frame. If the sum in the middle of (1)
is convergent in norm, the frame is called {\it standard}.

If $\{x_i\}_{i\in I}$ is a  standard frame in a finitely or
countably generated Hilbert $A$-module, it has a  unique operator
$S\in\ {\rm End}_A^*(H)$, where ${\rm End}_A^*(H)$ is the set of
adjointable $A$-linear maps on $H$, such that for every  $x\in H$,
\begin{equation*}
x=\sum_{i\in I}\langle  x,Sx_i\rangle  x_i=\sum_{i\in I}\langle
x,x_i\rangle  Sx_i.
\end{equation*}
Moreover $S$ is positive and invertible.}\vspace{.5pc}
\end{deff}

\begin{deff}$\left.\right.$\vspace{.5pc}

\noindent {\rm Let $H$ be a Hilbert $A$-module, and let $v\in H$.
We say that $v$ is a basic element if $e=\langle v,v\rangle$ is a
minimal projection in \,$A$, \,i.e. \,$eAe={\Bbb C}e$. A system
$\{v_\lambda\hbox{\rm :}\ \lambda\in \Lambda\}$ of basic elements
of}\pagebreak

\noindent {\rm $H$ is called orthonormal if $\langle
v_\lambda,v_\mu\rangle=0$ for all $\lambda\ne \mu$. An {\it
orthonormal basis} for $H$ is an orthonormal system which
generates a dense submodule of $H$.}
\end{deff}

\section{Main results}

Let $A$ and $B$ be $C^*$-algebras, $E$  a Hilbert $A$-module  and
let $F$ be a Hilbert $B$-module. We take $A\otimes B$ as the
completion of $A\otimes_{\rm alg} B$ with the spatial norm. Hence
$A\otimes B$ is a $C^*$-algebra and for every $a\in A$, $b\in B$
we have $\|a\otimes b\|=\|a\|\cdot\|b\|$.  The algebraic tensor
product $E\otimes_{\rm alg} F$ is a pre-Hilbert $A\otimes
B$-module with module action
\begin{equation*}
(a\otimes b)(x\otimes y)=ax\otimes by\qquad (a\in A,\ b\in B,\
x\in E,\ y\in F),
\end{equation*}
and $A\otimes B$-valued inner product
\begin{equation*}
\langle x_1\otimes y_1,x_2\otimes y_2\rangle  =\langle
x_1,x_2\rangle \otimes
 \langle  y_1,y_2\rangle \qquad (x_1,x_2\in E,\ y_1,y_2\in F).
 \end{equation*}
We also know that for $z=\sum_{i=1}^n x_i\otimes y_i$ in
$E\otimes_{\rm alg} F$ we have
\begin{equation*}
\langle z,z\rangle =\sum_{i,j}\langle  x_i,x_j\rangle  \otimes
 \langle  y_i,y_j\rangle  \geq 0
 \end{equation*}
and $\langle z,z\rangle  =0$ if and only if $z=0$. Just as in the
case of ordinary pre-Hilbert space, we can form the completion
$E\otimes F$ of $E\otimes_{\rm alg} F$, which is a Hilbert
$A\otimes B$-module. It is called the {\it tensor product} of $E$
and $F$ (see \cite{Lance}). We note that if $a\in A^+$ and $b\in
B^+$, then $a\otimes b\in (A\otimes B)^+$. Plainly if $a$, $b$ are
hermitian elements of $A$ and $a\geq b$, then for every positive
element $x$ of $B$, we have $a\otimes x\geq b\otimes x$.

\setcounter{defin}{0}
\begin{lem}
Let $\{u_i\}_{i\in I}$ be a frame for $E$ with frame bounds $A$
and $B${\rm ,} and let $\{v_j\}_{j\in J}$ be a frame for $F$ with
frame bounds $C$ and $D$. Then $\{u_i\otimes v_j\}_{i\in I, j\in
J}$ is a frame for $E\otimes F$ with frame bounds $AC$ and $BD$.
In particular{\rm ,} if $\{u_i\}_{i\in I}$ and $\{v_j\}_{j\in J}$
are tight or Parseval frames{\rm ,} then so is $\{u_i\otimes
v_j\}_{i\in I, j\in J}$.
\end{lem}

\begin{proof}
Let $ x\in E$ and $y \in F$. Then we have
\begin{align}
A\langle  x,x\rangle  &\leq  \sum_{i\in I}\langle  x,u_i\rangle
\langle
 u_i,x\rangle  \leq B\langle  x,x\rangle ,\\[.5pc]
C\langle  y,y\rangle  &\leq \sum_{j\in J}\langle  y,v_j\rangle
\langle v_j,y\rangle  \leq D\langle  y,y\rangle .
\end{align}
Therefore
\begin{align*}
A\langle  x,x\rangle  \otimes\langle  y,y\rangle  & \leq  \sum_i
\langle x,u_i\rangle  \langle  u_i,x\rangle  \otimes \langle  y,y\rangle\\[.5pc]
& \leq B\langle x,x\rangle  \otimes \langle  y,y\rangle  .
\end{align*}
Now by (3), we have
\begin{align*}
AC\langle  x,x\rangle  \otimes\langle  y,y\rangle  & \leq
\sum_i\sum_j \langle x,u_i\rangle  \langle  u_i,x\rangle  \otimes
\langle y,v_j\rangle  \langle v_j,y\rangle  \\[.5pc] & \leq   B\langle
x,x\rangle  \otimes\sum_j\langle
y,v_j\rangle  \langle  v_j,y\rangle  \\[.5pc]
& \leq  BD\langle  x,x\rangle  \otimes \langle  y,y\rangle  .
\end{align*}
Consequently we have
\begin{align*}
AC\langle  x\otimes y,x\otimes y\rangle  &\leq \sum_i\sum_j
\langle x\otimes y,u_i\otimes v_j\rangle  \langle  u_i\otimes
v_j,x\otimes y\rangle  \\[.5pc] &\leq BD\langle  x\otimes y,x\otimes
y\rangle .
\end{align*}
From these inequalities it follows that for all $z=\sum_{k=1}^n
x_k\otimes y_k$ in $E\otimes_{\rm alg} F$,
\begin{equation}
AC\langle  z,z\rangle    \leq  \sum_{i,j} \langle  z,u_i\otimes
v_j\rangle  \langle u_i\otimes v_j,z\rangle   \leq BD\langle
z,z\rangle  .
\end{equation}
Hence relation (4) holds for all $z$ in $E\otimes F$. \hfill
$\Box$
\end{proof}

From Theorem 1 of \cite{bak} and the above lemma we have the
following result.

\begin{theo}[\!]
Let $E$ be a Hilbert $A$-module and $F$ be a Hilbert $B$-module.
Let $\{u_i\}_{i\in I}$ and $\{v_j\}_{j\in J}$ be orthonormal bases
in $E$ and $F${\rm ,} respectively. Then $\{u_i\otimes v_j\}_{i\in
I, j\in J}$ is an orthonormal basis for $E\otimes F$.
\end{theo}

\begin{proof}
It is clear that each $u_i\otimes v_j$ is a basic element of
$E\otimes F$ and $\{u_i\otimes v_j\}_{i\in I, j\in J}$ is an
orthonormal system in $E\otimes F$. Now for each $x\in E$ and each
$y\in F$, we have $x=\sum_{i\in I}\langle x, u_i\rangle u_i$ and
$y=\sum_{j\in J}\langle y, v_j\rangle v_j$. Hence
\begin{equation*}
x\otimes y=\sum_{i\in I}\sum_{j\in J}\langle x\otimes y,
u_i\otimes v_j\rangle
  u_i\otimes v_j.
  \end{equation*}
Similar to the above lemma we can show that for
  each $z$ in $E\otimes F$, we have $z=\sum_{i\in I}\sum_{j\in
  J}\langle z, u_i\otimes v_j\rangle u_i\otimes v_j$. But Bakic and
  Guljas in Theorem 1 of \cite{bak} showed that if $W$ is a Hilbert
  $C^*$-module over a $C^*$-algebra $A$, and $(v_\lambda)_{\lambda\in \Lambda}$
  is an orthonormal system in $W$, then $(v_\lambda)_{\lambda\in
  \Lambda}$ is an orthonormal basis for $W$ if and only if for every
  $w\in W$, $w=\sum \langle w,v_\lambda\rangle v_\lambda$. Now  by using this
  fact we have the result. \hfill $\Box$
  \end{proof}

Let $\{u_i\}_{i\in I}$ and $\{v_j\}_{j\in J}$ be standard frames
for $E$ and $F$, respectively. So $\{u_i\otimes v_j\}_{i\in I,j\in
J}$ is a standard frame for $E\otimes F$.

Let $S$, $S'$ and $S''$ be the frame operators of $\{u_i\}_{i\in
I}$, $\{v_j\}_{j\in J}$ and $\{u_i\otimes v_j\}_{i\in I,j\in J}$,
respectively. So $S$ is $A$-linear and $S'$ is $B$-linear. Hence for
every $x\in E$ and $y\in F$, we have $x=\sum_i\langle  x,Su_i\rangle
u_i$, $y=\sum_j \langle  y,S'v_j\rangle  v_j$. Therefore
\begin{align*}
x\otimes y & =  \sum_i\sum_j \langle  x,Su_i\rangle  u_i \otimes
\langle
  y,S'v_j\rangle  v_j\\[.5pc]
& = \sum_i\sum_j (\langle  x,Su_i\rangle  \otimes \langle
y,S'v_j\rangle
 ) (u_i\otimes v_j)\\[.5pc]
& = \sum_i\sum_j (\langle  x\otimes y, Su_i\otimes S'v_j\rangle
u_i\otimes v_j.
\end{align*}
Now by the uniqueness of frame operator we have $S''(u_i\otimes
v_j)=Su_i\otimes S'v_j$. Hence $S''=S\otimes S'$, which is a
bounded $A\otimes B$-linear, self-adjoint, positive
 and invertible operator on $E\otimes F$.  We note that
 $\|S''\|=\|S\otimes S'\|\leq \|S\|.\|S'\|$.
Now we summarize the above results as follows:

\begin{theo}[\!]
Let $\{u_i\}_{i\in I}$  and  $\{v_j\}_{j\in J}$   be standard
frames in the Hilbert $C^*$-modules $E$ and $F${\rm ,}
respectively. If $S${\rm ,} $S'$ and $S''$ are the frame operators
of $\{u_i\}_{i\in I}${\rm ,}  $\{v_j\}_{j\in J}$ and $\{u_i\otimes
v_j\}_{i\in I,j\in J}${\rm ,} respectively{\rm ,} then
$S''=S\otimes S'$. \end{theo}

For the frame operator we prove the following result.

\begin{lem}
If $\{x_i\}_{i\in I}$ is a frame in Hilbert $A$-module $X$ with
frame operator $S$ and $Q\in\ {\rm End}_A^*(X)$  is invertible{\rm
,} then $\{Qx_i\}_{i\in I}$ is a frame  in $X$ with frame operator
$Q^{*-1}SQ^{-1}$.
\end{lem}

\begin{proof}
Let $\{x_i\}_{i\in I}$ be a frame of $X$ with frame
 operator $S$. Then there exist constants $A$, $B>0$ such that for
 every $x\in X$,
\begin{equation}
 A\langle x,x\rangle\leq \sum_i|\langle x,x_i\rangle|^2\leq
 B\langle x,x\rangle,
\end{equation}
and $S^{-1}x=\sum_i\langle  x,x_i\rangle  x_i$. Since $Q$ is
invertible and  $Q\in\ \hbox{End}^*_A(X)$, then $Q$ is  a bounded
$A$-linear map with  invertible adjoint $Q^*$. So for every $x\in
X$, we have
\begin{equation}
\|Q^{*-1}\|^{-1}\cdot|x|\leq |Q^*x|\leq \|Q^*\|\cdot |x|.
\end{equation}
Since $Q$ is $A$-linear, $QS^{-1}x=\sum_i\langle  x,x_i\rangle
Qx_i$. So $QS^{-1}Q^*(Q^{*-1}x)=\sum_i\langle
Q^{*-1}x,Qx_i\rangle$ $Qx_i$,  because
\begin{equation*}
\langle  x,x_i\rangle  =\langle
Q^*Q^{*-1}x,x_i\rangle  =\langle
 Q^{*-1}x,Qx_i\rangle.
 \end{equation*}
Consequently, for every $x\in X$,
\begin{equation}
QS^{-1}Q^*(x)=\sum_i\langle x,Qx_i\rangle  Qx_i.
\end{equation}
Now by using (5) and (6) we have
\begin{align*}
A\|Q^{*-1}\|^{-2} \langle x,x\rangle  &\leq A\langle
Q^*x,Q^*x\rangle\\[.5pc]
&\leq  \sum_i|\langle Q^*x,x_i\rangle  |^2  \leq B\langle
Q^*x,Q^*x\rangle \leq B\|Q^*\|^2\langle x,x \rangle .
\end{align*}
On the other hand, $\langle Q^*x,x_i\rangle=\langle
x,Qx_i\rangle$, so $\{Qx_i\}_{i\in I}$ is a frame for $X$ and by
(7), $Q^{*-1}SQ^{-1}=(QS^{-1}Q^*)^{-1}$ is the  frame operator of
$\{Qx_i\}_{i\in I}$. \hfill $\Box$
\end{proof}

\begin{theo}[\!]
If $Q\in {\rm End}_A^*(E)$ is an invertible $A$-linear map and
$\{T_i\}_{i\in J}$ is a frame in $E\otimes F$ with frame operator
$S${\rm ,} then $\{(Q^*\otimes I)(T_i)\}_{i\in J}$ is a frame of
$E\otimes F$ with frame operator $(Q\otimes I)^{-1}S(Q^*\otimes
I)^{-1}$.
\end{theo}

\begin{proof}
Since $Q\in \hbox{End}^*_A(E)$, $Q\otimes I\in
\hbox{End}_A^*(E\otimes F)$ with inverse $Q^{-1}\otimes I$. It is
obvious that $Q\otimes I$ is $A\otimes B$-linear, adjointable,
with adjoint $Q^*\otimes I$. An easy calculation shows that for
every elementary tensor $x\otimes y$,
\begin{align*}
\|(Q\otimes I)(x\otimes y)\|^2 & =  \|Q(x)\otimes
y\|^2=\|Q(x)\|^2\cdot \|y\|^2\\[.5pc]  & \leq  \|Q\|^2\cdot \|x\|^2\cdot
\|y\|^2=\|Q\|^2\cdot \|x\otimes y\|^2.
\end{align*}
So $Q\otimes I$ is bounded, and therefore it can be extended to
$E\otimes F$. Similarly for $Q^{*-1}\otimes I$. Hence $Q\otimes I$
is $A\otimes B$-linear, adjointable with adjoint $Q^*\otimes I$,
and as we mentioned in the proof of Lemma 3.4, $Q^*$ is invertible
and bounded. Hence for every $T\in E\otimes F$, we\break have
\begin{equation}
\|Q^{*-1}\|^{-1}\cdot|T|\leq |(Q^*\otimes I)T|\leq \|Q\|\cdot |T|.
\end{equation}
Hence $Q\otimes I\in \hbox{End}^*_{A\otimes B}(E\otimes F)$. Now
by the above lemma we have the result.\hfill $\Box$
\end{proof}

Now we generalize some of the results in \cite{ka} to frame of
subspaces. First we recall the definition of frame of subspaces
(for basic definitions  and properties, see \cite{cas.sub}).

\begin{deff}$\left.\right.$\vspace{.5pc}

\noindent {\rm Let $H$ be a separable Hilbert  space and let
$\{v_i\}_{i\in I}$ be a  sequence of weights, i.e., $v_i>o$ for
all $i \in I$. A sequence $\{W_i\}_{i\in I}$ of closed subspaces
of $H$ is a {\it frame of  subspaces} with respect  to
$\{v_i\}_{i\in I}$ if there exist real numbers $A,B>0$ such that
for every $x\in H$,
\begin{equation*}
A\|x\|^2\leq \sum_{i\in I} v_i^2\|\pi_{W_i}(x)\|^2\leq B\|x\|^2,
\end{equation*}
where for each $i\in I$, $\pi_{W_i}$ is the orthogonal projection
of $H$ onto $W_i$. Similar to frames, $A$ and $B$ are called the
frame bounds. 1f $A=B=\lambda$, the  frame of subspaces is
$\lambda$-tight and it is a Parseval frame of subspaces if
$A=B=1$.

Let $H$ and $K$ be Hilbert spaces and let $W$, $Z$ be closed
subspaces of $H$ and $K$, respectively. Then $\pi_W\otimes
\pi_Z\hbox{\rm :}\ H\otimes_{\rm alg}K\rightarrow W\otimes Z$ is a
bounded linear map, and it can be extended to a bounded linear map
from $H\otimes K$ into $W\otimes Z$.  We also denote it by
$\pi_W\otimes \pi_Z$ and clearly it is surjective. Hence
$\pi_W\otimes \pi_Z$  is the orthogonal projection of $H\otimes K$
onto $W\otimes Z$.}
\end{deff}

\begin{theo}[\!]
Let $\{W_i\}_{i\in I}$ be a frame of subspaces with respect to
$\{u_i\}_{i\in I}$  for $H${\rm ,} with frame bounds $A${\rm ,}
$B${\rm ,} and let $\{Z_j\}_{j\in J}$ be a frame of subspaces with
respect to $\{v_j\}_{j\in J}$  for $K$ with frame bounds $A'${\rm
,} $B'$. Then $\{W_i\otimes Z_j\}_{i\in I,j\in J}$ is a frame of
subspaces with respect to $\{u_iv_j\}_{i\in I, j\in J}$ for
$H\otimes K$ with frame bounds $AA'$ and $BB'$. It is tight or
Parseval if $\{W_i\}_i$ and  $\{Z_j\}$ are tight or Parseval.
\end{theo}

\begin{proof}
Let $x\otimes y$  be an elementary tensor. Then $A\|x\|^2\leq
\sum_{i\in I} u_i^2\|\pi_{W_i}(x)\|^2\leq B\|x\|^2$ and
$A'\|y\|^2\leq \sum_{j\in J} v_j^2\|\pi_{Z_j}(y)\|^2\leq
B'\|y\|^2$.

A simple calculation shows that
\begin{align*}
AA'\|x\otimes y\|^2& \leq  \sum_i\sum_j
u_i^2v_j^2\|\pi_{W_i}(x)\|^2\cdot \|\pi_{Z_j}(y\|^2\\[.5pc]
& \leq  BB'\|x\otimes y\|^2.
\end{align*}
Hence
\begin{equation*}
AA'\|x\otimes y\|^2\leq \sum_{i,j}u_i^2v_j^2\|\pi_{W_i}(x)\otimes
\pi_{Z_j}(y)\|^2\leq BB'\|x\otimes y\|^2.
\end{equation*}
Therefore
\begin{equation}
AA'\|x\otimes y\|^2\leq \sum_{i,j}u_i^2v_j^2\|\pi_{W_i}\otimes
\pi_{Z_j}(x\otimes y)\|^2\leq BB'\|x\otimes y\|^2.
\end{equation}
Consequently, for every $z=\sum_{l=1}^n x_l\otimes y_l$ in
$H\otimes_{\rm alg} K$  and every $z$ in $H\otimes K$, the
relation (9) holds. Hence we have the result. \hfill $\Box$
\end{proof}

Now we try to generalize a known result of frames (Proposition~3.1
of \cite{ka}) to frames of subspaces.

\begin{deff}$\left.\right.$\vspace{.5pc}

\noindent {\rm Let $\{W_i\}_{i\in I}$ be a  frame of subspaces for
$H$ with respect to $\{v_i\}_{i\in I}$. Then the frame operator
$S_{W,v}$ for $\{W_i\}_{i\in I}$ and $\{v_i\}_{i\in I}$ is defined
by
\begin{equation*}
S_{W,v}(x)=\sum_{i\in I}v_i^2\pi_{W_i}(x),\quad x\in H
\end{equation*}}
\end{deff}

\begin{coro}$\left.\right.$\vspace{.5pc}

\noindent With the hypothesis in Theorem {\rm 3.7,} if $S_{W,u}$
and $S_{Z,v}$ are frame operators for $\{W_i\}_{i\in I}${\rm ,}
 $\{u_i\}$ and $\{Z_j\}${\rm ,} $\{v_j\}${\rm ,} respectively{\rm ,} then $S_{W,u}\otimes S_{Z,v}$
is the frame operator for $\{W_i\otimes Z_j\}_{i\in I,j\in J}$ and
$\{u_iv_j\}_{i\in I,j\in J}$.
\end{coro}

\begin{proof}
Let $x\otimes y$ be an elementary tensor. Therefore
\begin{align*}
S_{W,u}\otimes
S_{Z,v}(x\otimes y)& = S_{W,u}(x)\otimes S_{Z,v}(y)\\[.5pc]
& = \sum_iu_i^2\pi_{W_i}(x)\otimes
\sum_jv_j^2\pi_{Z_j}(y)\\[.5pc]
&=\sum_i\sum_ju_i^2v_j^2(\pi_{W_i}\otimes \pi_{Z_j})(x\otimes y).
\end{align*}
Now the uniqueness of frame operator implies that $S_{W,u}\otimes
S_{Z,v}$ is the desired frame operator.\hfill $\Box$
\end{proof}

\begin{rem}
{\rm Let $H$ and $K$ be Hilbert spaces. A map $T\hbox{\rm :}\
H\longrightarrow K$ is {\it antilinear} (or conjugate linear) if
$T(\lambda x+y)=\bar{\lambda}T(x)+T(y)$ for all $\lambda \in {\Bbb
C}$ and $x,y\in H$. By  the techniques in \cite{folland},
$H\otimes K$ is the set of anti-linear maps $T\hbox{\rm :}\
K\rightarrow H$ with the norm $\|.\|$ defined by
\begin{equation*}
\|T\|=\hbox{sup}\{\|Ty\|\hbox{\rm :}\ y\in K,\ \|y\|\leq 1\}.
\end{equation*}
So $W_i\otimes Z_j$ is the  set of anti-linear maps $T\hbox{\rm
:}\ Z_j\rightarrow W_i$ and  therefore $\pi_{W_i}\otimes
\pi_{Z_j}$
 is the map which assigns to every $T\in H\otimes K$, the restriction of
$\pi_{W_i}\circ T$ to $Z_j$, i.e. $\pi_{W_i}\circ T|Z_j$.}
\end{rem}

\section{Resolution of the identity}

In this section we present the notion of $\ell^2$-resolution of
the identity with lower resolution bound in tensor product of
Hilbert spaces (for more information see
\cite{cas.sub,fornasier}).

\setcounter{defin}{0}
\begin{deff}$\left.\right.$\vspace{.5pc}

\noindent {\rm Let $I$ be a countable index set and let $H$ be a
Hilbert space. Let $\{v_i\}_{i\in I}$ be a family of weights,
i.e., for all $i$, $v_i>0$. Then a family of bounded operators
$\{T_i\}_{i\in I}$ on $H$ is called a {\it $\ell^2$-resolution of
the identity with lower resolution bound} with respect to
$\{v_i\}_{i\in I}$ on  $H$ if there are positive real numbers $C$
and $D$ such that for all $f\in H$,

\begin{enumerate}
\renewcommand{\labelenumi}{(\roman{enumi})}
\leftskip .1pc
\item $C\|f\|^2\leq \sum_{i\in I}v_i^{-2}\|T_i(f)\|^2\leq
D\|f\|^2$,

\item $f=\sum_{i\in I}T_i(f)$ (and the series converges
unconditionally for every $f\in H$).\vspace{-.3pc}
\end{enumerate}

The optimal values of $C$ and $D$ are called the {\it bounds} of
the resolution of the identity.}
\end{deff}

\begin{propo}$\left.\right.$\vspace{.5pc}

\noindent Let $\{T_i\}_{i\in I}$ be a $\ell^2$-resolution of the
identity with lower resolution bound with  respect to
$\{v_i\}_{i\in I}$ on $H${\rm ,} and let $\{S_j\}_{j\in J}$ be a
$\ell^2$-resolution of the identity with lower resolution bound
with respect to $\{u_j\}_{j\in J}$ on $K$. Then $\{T_i\otimes
S_j\}_{i\in I,j\in J}$ is a $\ell^2$-resolution of the identity
with lower resolution bound with respect to $\{v_iu_j\}_{i\in
I,j\in J}$ on $H\otimes K$.
\end{propo}

\begin{proof}
Let $f\in H$, $g\in K$. Then $f=\sum_{i\in I}T_i(f)$,
$g=\sum_{j\in J}S_j(g)$, and consequently
\begin{align*}
\sum_{i,j}(T_i\otimes S_j)(f\otimes g)& = \sum_{i,j}T_i(f)\otimes
S_j(g)\\[.5pc]
& =  \sum_iT_i(f)\otimes \sum_jS_j(g)=f\otimes g.
\end{align*}
Since both the series  $f=\sum_{i\in I}T_i(f)$ and $g=\sum_{j\in
J}S_j(g)$  are  unconditionally convergent, the above series is
unconditionally convergent. So for every $h\in H\otimes_{\rm alg}
K$ and consequently for every $h\in H\otimes K$ the above relation
holds. Let $C$, $D$ and $C'$, $D'$ be the bounds of the
resolutions $\{T_i\}$ and $\{S_j\}$,
  respectively. Then for every $f\in H$, $g\in K$  we have
\begin{align}
CC'\|f\otimes g\|^2 & = CC'\|f\|^2\cdot \|g\|^2\leq
C'\sum_iv_i^{-2}\|T_i(f)\|^2\cdot \|g\|^2\nonumber\\[.5pc]
& \leq \sum_i v_i^{-2}\|T_i(f)\|^2\cdot \sum_ju_j^{-2}\|S_jg\|^2 \nonumber\\[.5pc]
&= \sum_{i,j} v_i^{-2}u_j^{-2}\|(T_i\otimes S_j)(f\otimes g)\|^2 \nonumber\\[.5pc]
& \leq  DD'\|f\otimes g\|^2.
\end{align}
Now by using the fact that
\begin{equation*}
\|(T\otimes S)\left(\sum_{i=1}^n f_i\otimes
g_i\right)\|^2=\|T\left(\sum_{i=1}^n
f_i\right)\|^2\cdot\|S\left(\sum_{i=1}^n g_i\right)\|^2,
\end{equation*}
and $\|\sum_{i=1}^n f_i\otimes g_i\|^2=\|\sum_{i=1}^n f_i\|^2\cdot
\|\sum_{i=1}^n g_i\|^2$, we conclude that  for every
$h=\sum_{i=1}^n f_i\otimes g_i$ and consequently for every $h\in
H\otimes K$ the relation (10) holds. \hfill $\Box$
\end{proof}

From  the above proposition and Proposition 3.26 of \cite{cas.sub}
we have the following result.

\begin{coro}$\left.\right.$\vspace{.5pc}

\noindent With the hypothesis in Corollary {\rm 3.9,} if
$T_i=\pi_{W_i}S_{W,v_i}$ and $S_j=\pi_{Z_j}S_{Z,u_j}${\rm ,} then
$\{v_i^2u_j^2T_i\otimes S_j\}_{i\in I,j\in J}$ is a
$\ell^2$-resolution of the identity with lower resolution bound
with respect to $\{v_iu_j\}_{i\in I, j\in J}$ on $H\otimes K$ and
for all $z\in H\otimes K\!{\rm ,}$\pagebreak
\begin{equation*}
\frac{C}{D^2}\cdot \frac{C'}{D'{^2}}\|z\|^2\leq \sum_{i\in
I}\sum_{j\in J}v_i^2u_j^2\|(T_i\otimes S_j)(z)\|^2\leq
\frac{D}{C^2}\cdot \frac{D'}{C'{^2}}\|z\|^2.
\end{equation*}
\end{coro}

\section{Frame representation}

Let $H$ be a separable Hilbert space, and let $G$ be a discrete
countable abelian group. Let $\pi\hbox{\rm :}\ G\rightarrow B(H)$
be a unitary representation of $G$ on $H$. If there is a vector
$v\in H$ such that $\{\pi(g)v|g\in G\}$ is a frame for $H$, then
the representation $\pi$ is called a {\it frame representation}.
Let $\hat{G}$ denote the dual group of $G$, i.e., the group of
characters on $G$ and let $\lambda$ be the normalized Haar measure
on $\hat{G}$. Let $\pi:G\longrightarrow B(H)$  be a frame
representation with frame vector $v$. As we have in
\cite{altw,HL,mur} there is a spectral measure $E$ on $\hat{G}$
such that
\begin{equation*}
\pi(g)=\int_{\hat{G}} g(\xi)\hbox{d}E(\xi).
\end{equation*}
Since $\pi$ is a frame representation, by using the results in \S
2 of \cite{altw} and  the properties of spectral measure there is
a unitary operator $U\hbox{\rm :}\ H\longrightarrow L^2(F,
\lambda|F)$, where $F$ is a measurable subset of $\hat{G}$ with
$\lambda(F)>0$ and $\lambda|F$ is the restriction of Haar measure
$\lambda $ to $F$ such that $U$ interwines the spectral measure on
$H$ and the canonical spectral measure on $\hat{G}$. The operator
$U$ is called the {\it decomposition operator}.
 Moreover $\pi$ is unitarily equivalent to the
representation $\sigma\hbox{\rm :}\ G\longrightarrow
B(L^2(F,\lambda|F))$ defined by $\sigma(g)=M_g$, where $M_g$ is
the multiplication operator with symbol $g$. In fact,
$U^*M_gU=\pi(g)$.

We also note that if  $\theta_v$ is the analysis operator of $H$
for frame vector $v$, then $\theta_v\pi(g)=L_g\theta_v$, where
$L_g\hbox{\rm :}\ \ell^2(G)\longrightarrow \ell^2(G)$ is defined
by $(L_gx)(h)=x(g^{-1}h)$ for all $h\in G$. In fact, if $J$ is the
range of $\theta_v$, then the representation $\pi$ of $G$ is
unitarily equivalent to $\rho=L_g|J$ (see Lemma 3 of \cite{altw}).
For more details see \cite{altw} or \cite{HL}.

Let  $H$ and $K$ be separable Hilbert spaces and let $\pi\hbox{\rm
:}\ G_1\rightarrow B(H)$ and $\sigma\hbox{\rm :}\ G_2\rightarrow
B(K)$ be frame representations on $H$ and $K$ with frame vectors
$v\in H$ and $w\in K$, respectively. Since $G_1$ and $G_2$ are
discrete countable abelian groups, their direct sum $G=G_1\oplus
G_2$ is a discrete countable abelian group. Hence we can consider
the representation $\pi \otimes \sigma\hbox{\rm :}\ G\rightarrow
B(H\otimes K)$ defined by
\begin{equation*}
(\pi\otimes \sigma)(g,h)=\pi g\otimes \sigma h,\quad (g,h)\in G.
\end{equation*}
Since $\{\pi(g)v\hbox{\rm :}\ g\in G_1\}$  is a frame for $H$ and
$\{\sigma(h)w\hbox{\rm :}\ h\in G_2\}$ is a frame for $K$, by
Lemma 3.1 and the definition of $\pi\otimes \sigma$,
\begin{equation*}
\{\pi\otimes\sigma(g,h)(v\otimes w): (g,h)\in G\}= \{(\pi
g)v\otimes(\sigma h)w: (g,h)\in G\}
\end{equation*}
is a frame for $H\otimes K$. So $\pi\otimes \sigma$ is a frame
representation of $H\otimes K$ with frame vector~$v\otimes w$.
Moreover, if $\theta_v$ and $\theta_w$ are the analysis operators
of $H$ and $K$ for frame vectors $v$ and $w$, respectively, then
$\theta_v\otimes \theta_w$ is the analysis operator of $H\otimes
K$
 for frame vector $v\otimes w$. Hence we have proved the
following result.

\setcounter{defin}{0}
\begin{theo}[\!]
Let $\pi\hbox{\rm :}\ G_1\rightarrow  B(H)$ and $\sigma\hbox{\rm
:}\ G_2\rightarrow B(K)$
 be frame representations with frame vectors $v$ and $w${\rm ,} respectively.
Then $\pi\otimes \sigma\hbox{\rm :}\ G_1\oplus G_2\rightarrow
B(H\otimes K)$ is a frame representation with frame vector
$v\otimes w$. If $\theta_v$ and $\theta_w$  are the analysis
operators for frame vectors $v$ and $w${\rm ,} respectively{\rm ,}
then $\theta_v\otimes \theta_w$ is the analysis operator for
$v\otimes w$. \hfill $\Box$
\end{theo}

For the decomposition operators we have the following result.

\begin{theo}[\!]
With the hypothesis in Theorem {\rm 5.1,} suppose that $U\hbox{\rm
:}\ H\rightarrow L^2(E,\lambda|E)$ and $V\hbox{\rm :}\
K\rightarrow L^2(F,\lambda|F)$ are the decomposition operators of
$\pi$ and $\sigma${\rm ,} respectively{\rm ,} then $U\otimes
V\hbox{\rm :}\ H\otimes K\rightarrow L^2(E\oplus F,\lambda\times
\mu |E\otimes F)$ is the decomposition operator of $\pi\otimes
\sigma$.
\end{theo}

\begin{proof}
It is clear that $(G_1\oplus G_2)^\wedge=\hat{G_1}\oplus
\hat{G_2}$. If $U\hbox{\rm :}\ H\longrightarrow L^2(E,\lambda |E)$
and $V\hbox{\rm :}\ K\longrightarrow L^2(F,\mu |F)$, where
$\hat{G_1}\supseteq E$, $\hat{G_2}\supseteq F$, then
$\hat{G_1}\oplus \hat{G_2}\supseteq E\oplus F$ and $U\otimes
V:H\otimes K\longrightarrow L^2(E\oplus F,\lambda\times \mu
|E\oplus F)$, where $\lambda\times \mu$  is the product measure of
$\lambda$ and $\mu$.  We note that for every $x\in H$, $y\in K$,
the function $(U\otimes V)(x\otimes y)=Ux\otimes Vy$ defined  on
$E\oplus F$ by $(Ux\otimes Vy)(\zeta,\eta)=(Ux)(\zeta).(Vy)(\eta)$
and since $L^2(E,\lambda|E)\otimes L^2(F,\lambda|F)$ is isomorphic
to $L^2(E\oplus F,\lambda\times\mu|E\oplus F)$ we can  take
$Ux\otimes Vy\in L^2(E\oplus F,\lambda\times\mu|E\oplus F)$. Since
$G_1$ and $G_2$ form an orthonormal basis of
$L^2(\hat{G}_1,\lambda)$ and $L^2(\hat{G}_2,\mu)$, respectively
(Corollary 4.26 of \cite{folland}), a simple calculation shows
that
\begin{align*}
\|Uv\otimes Vw\|^2 &=\|\chi_{E\oplus F}\cdot Uv\otimes Vw\|^2\\[.5pc]
& =\int_{\hat{G}_1}|\chi_E(\zeta)Uv(\zeta)|^2\hbox{d}\lambda\cdot
\int_{\hat{G}_2}|\chi_F(\eta)Vw(\eta)|^2\hbox{d}\mu\\[.5pc]
&=\|\chi_EUv\|^2\cdot \|\chi_FVw\|^2=\|Uv\|^2\cdot \|Vw\|^2<\infty.\\[-3pc]
\end{align*}
\hfill $\Box$
\end{proof}

\begin{coro}$\left.\right.$\vspace{.5pc}

\noindent Let $\{\pi(g)v\}_{g\in G_1}$ and $\{\sigma(h)w\}_{h\in
G_2}$  be frames for $H$ and $K$ with frame bounds  $A_1$, $B_1$
and $A_2${\rm ,} $B_2${\rm ,} respectively. Then $\{(\pi\otimes
\sigma)(g,h)(v\otimes w)\}_{g\in G_1,h\in G_2}$ is a frame with
frame bounds $A_1A_2$ and $B_1B_2$.
\end{coro}

\begin{proof}
First we note that for all $x\in H$,
\begin{equation*}
\sum_{g\in G_1} |\langle x,\pi(g)v\rangle|^2=\sum_{g\in G}
\int_{\hat{G}_1}
|Ux(\zeta)Uv(\zeta)|^2\hbox{d}\lambda=\|(Ux)(Uv)\|^2
\end{equation*}
and
\begin{equation*}
A_1\|x\|^2\leq \sum_{g\in G_1}|\langle x,\pi(g)v\rangle|^2\leq
B_1\|x\|^2, \ \ \mbox{ for all } x\in H.
\end{equation*}
Similarly
\begin{equation*}
A_2\|y\|^2\leq \sum_{h\in G_2}|\langle y,\sigma(h)w\rangle|^2\leq
B_2\|y\|^2, \ \ \mbox{ for all } y\in K.
\end{equation*}
Hence for every elementary tensor $x\otimes y$ we have $\|x\otimes
y\|=\|x\|.\|y\|$ and
\begin{align*}
&\sum_{g\in G_1}\sum_{h\in G_2}|\langle x\otimes y,\pi(g)\otimes
\sigma(h)(v\otimes w)\rangle|^2\\[.5pc]
&\quad\,=\int_{\hat{G}_1}\int_{\hat{G}_2}|v(x)|^2\cdot |Uv|^2\cdot
|w(x)|^2\cdot
|Vw|^2d(\lambda\times \mu)\\[.5pc]
&\quad\,=\|(Ux)(Uv)\|^2\cdot \|(Vy)(Vw)\|^2.
\end{align*}
So we have the result. \hfill $\Box$
\end{proof}

We can also state similar results for Bessel vectors.

\begin{deff}$\left.\right.$\vspace{.5pc}

\noindent {\rm Let $\pi\hbox{\rm :}\ G\longrightarrow B(H)$  be a
frame representation with frame vector $v$. We say $v'\in H$ is a
Bessel vector for the frame representation if there exists $C_2>0$
such that for all $x\in H$,
\begin{equation*}
\sum_{g\in G}|\langle x,\pi(g)v'\rangle|^2\leq C_2\|x\|^2.
\end{equation*}}
\end{deff}

\begin{lem}
Suppose $\pi$  and $\sigma$ are frame representations on $H$ and
$K$ with frame vectors $v$ and $w${\rm ,} respectively. If $v'$
and $w'$ are Bessel vectors for $\pi$ and $\sigma${\rm ,}
respectively{\rm ,} then $v'\otimes w'$ is a  Bessel vector for
$\pi\otimes \sigma$.
\end{lem}

\begin{proof}
By Theorem 5.1, $\pi\otimes \sigma$ is a frame
 representation with frame vector $v\otimes w$, and
 since  $v'$ and $w'$ are Bessel
vectors for $\pi$  and $\sigma$, respectively, there are  constants
$C_2$ and $C'_2$ such that
\begin{align*}
\sum_{g\in \hat{G}_1} |\langle x,\pi(g)v'\rangle|^2 &\leq
C_2\|x\|^2,\quad x\in H,\\[.5pc]
\sum_{h\in \hat{G}_2} |\langle y,\sigma(h)w'\rangle|^2 &\leq
C'_2\|y\|^2,\quad y\in K.
\end{align*}
Hence for every elementary tensor $x\otimes y$ we  have
\begin{equation*}
\sum_{g\in \hat{G}_1}\sum_{h\in \hat{G}_2} |\langle x\otimes
y,\pi\otimes \sigma(g,h)(v'\otimes w')\rangle|^2\leq
C_2C'_2\|x\otimes y\|^2.
\end{equation*}
As we have in \$ 4, the above relation holds for every
$z=\sum_{i=1}^n x_i\otimes y_i$ and so for every $z\in H\otimes
K$. Therefore $ v'\otimes w'$ is a Bessel vector for $\pi\otimes
\sigma$. \hfill $\Box$
\end{proof}

\section*{Acknowledgement}

The authors express their gratitude to the referee for carefully
reading and  several valuable pointers which improved the
manuscript.

\end{document}